\documentclass{article}
\usepackage{amssymb}
\usepackage{amsmath}
\usepackage{amsthm}
\usepackage{mathtools}
\usepackage{enumerate}
\usepackage{tikz}
\usetikzlibrary{calc}
\usepackage{float}

\newtheorem*{thmHK}{Theorem (Heckman and Krakovski)}
\newtheorem*{problem}{Problem}

\newcommand{\EGC}{Erd\H{o}s-Gy\'{a}rf\'{a}s Conjecture}

\begin{document}

\title{Three Graphs and the \EGC}

\author{
Geoffrey Exoo \\
Department of Mathematics and Computer Science \\
Indiana State University \\
Terre Haute, IN 47809 \\
ge@cs.indstate.edu \\
}

\date{Dec 24, 2013}

\maketitle

\begin{abstract}
Three graphs related to the \EGC\, are presented.
The graphs are derived from the Buckyball, the Petersen graph,
and the Tutte-Coxeter graph.
The first graph is a partial answer to a question posed
by Heckman and Krakovski \cite{planar} in their recent work on the planar
version of the conjecture.  The other two graphs appear to be the
smallest known cubic graphs with no $2^m$-cycles for $m \leq 4$ and
for $m \leq 5$.
\end{abstract}

\section{\bf Introduction.}

The \EGC\, asserts that every graph with minimum degree at least $3$ contains
a cycle whose length is a power of $2$.
Recently, Heckman and Krakovski \cite{planar} proved the conjecture for
$3$-connected cubic planar graphs.

\begin{thmHK}
Every $3$-connected planar graph contains a $2^m$-cycle, for some $m \leq 7$.
\end{thmHK}

The authors suggest that the upper bound of $7$ might be improved, and
that a value as low as $4$ may be possible.
In this note, we show that the bound on $m$ must be at least $5$
by constructing a $3$-connected
cubic planar graph with neither $4$, $8$, nor $16$ cycles.
Then we consider the general problem of finding, for a given integer $k$,
the smallest $2^m$-cycle
free cubic graphs for all $m < k$, and
describe graphs for $k=4$ and $k=5$ that
appear to be the smallest known examples.

\section{\bf A $3$-connected cubic planar graph with no cycles of length $4$, $8$ or $16$.}

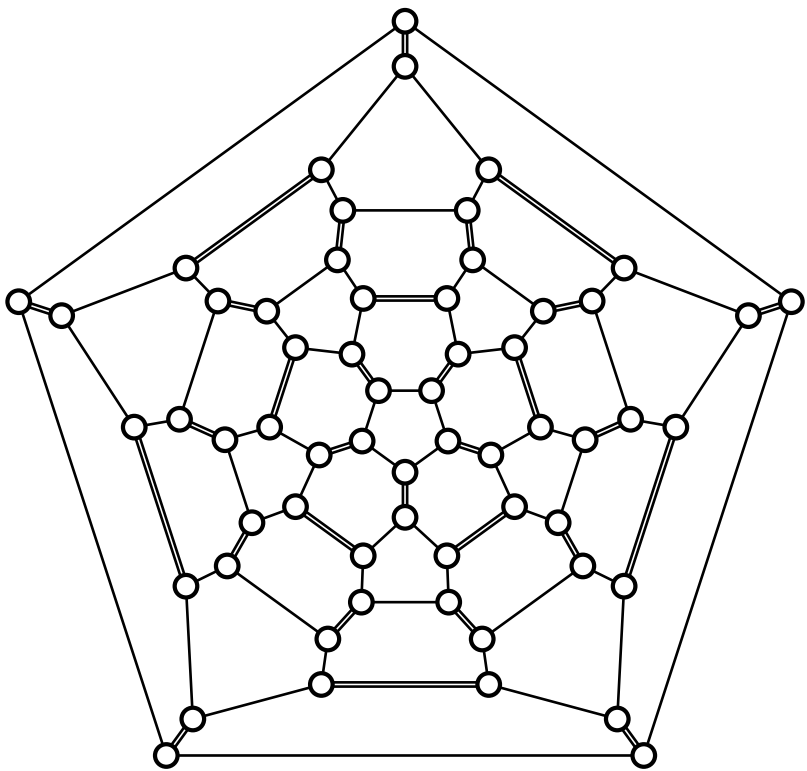
\begin{figure}[H]
\centering
\begin{tikzpicture}[line width=1pt,scale=0.60]
\tikzstyle{every node}=[draw=black,fill=white,ultra thick,
  shape=circle,minimum height=0.3cm,inner sep=2];

\foreach \theta in {0,72,...,288}
{
  \begin{scope}[rotate=\theta]
    \draw (54:1) -- (126:1);
    \draw[style=double] (72:3) -- (108:3);
    \draw (74:5) -- (106:5);
    \draw (18:9) -- (90:9);
    \draw[style=double] (54:1) -- (54:2);
    \draw (18:8) -- (36:6) -- (34:5);
    \draw[style=double] (34:5) -- (40:4);
    \draw (40:4) -- (36:3) -- ( 54:2) -- (72:3) -- (68:4);
    \draw[style=double] (68:4) -- (74:5);
    \draw (74:5) -- (72:6) -- (90:8);
    \draw[style=double] (90:8) -- (90:9);
    \draw[style=double] (36:6) -- (72:6);
    \draw (40:4) -- (68:4);
  \end{scope}
}

\foreach \theta in {0,72,...,288}
{
  \begin{scope}[rotate=\theta]
    \node[fill=white] at (54:1) {};
    \node[fill=white] at (54:2) {};
    \node[fill=white] at (72:3) {};
    \node[fill=white] at (36:3) {};
    \node[fill=white] at (68:4) {};
    \node[fill=white] at (40:4) {};
    \node[fill=white] at (74:5) {};
    \node[fill=white] at (34:5) {};
    \node[fill=white] at (72:6) {};
    \node[fill=white] at (36:6) {};
    \node[fill=white] at (90:8) {};
    \node[fill=white] at (90:9) {};
  \end{scope}
}

\end{tikzpicture}
\caption{The graph of the buckyball.}
\label{buckyball}
\end{figure}

The construction is based on the truncated icosahedron, known to chemists as $C_{60}$,
the buckminsterfullerene \cite{bucky}, or more popularly as the {\em buckyball}.
We recall some of the well-known graph theoretic properties of $C_{60}$.
It is a $3$-connected cubic planar vertex-transtive graph of order $60$.
It contains twelve $5$-cycles and twenty $6$-cycles, all of
which border faces in a plane drawing of the graph, such as
the one in Figure~\ref{buckyball}.

Under the action of the automorphism group there are two edge orbits.
One orbit contains $60$ edges, each of which
borders one pentagonal face and one hexagonal face.
These edges are called {\em single-bond} edges by chemists, a terminology
we also adopt.
The other edge orbit contains the $30$ edges that border two hexagonal
faces.  These are called {\em double-bond} edges.
Each vertex is incident with two single-bond edges and one
double-bond edge.
It will also be important to note that $C_{60}$ contains no other cycle
of length less than nine.

\begin{figure}[H]
\centering
\begin{tikzpicture}[line width=2pt,scale=2.0]
\tikzstyle{every node}=[draw=black,fill=white,ultra thick,
  shape=circle,minimum height=0.4cm,inner sep=2];
\def\xa{{1.0 - cos(72)}}
\def\xb{{2.0 + cos(72)}}
\def\ya{{sin(72)}}
\def\yb{{cos(36) + 1.0/(2.0 * cos(54))}}

  \draw (\xa,\ya) -- (0,0) -- (3,0) -- (\xb,\ya);
  \draw (1,0) -- (\xa,\ya) -- (1.5,\yb) -- (\xb,\ya) -- (2,0);

  \node[fill=white] at (0,0) {v};
  \node[fill=white] at (1,0) {};
  \node[fill=white] at (2,0) {};
  \node[fill=white] at (3,0) {w};
  \node[fill=white] at (\xa,\ya) {};
  \node[fill=white] at (1.5,\yb) {u};
  \node[fill=white] at (\xb,\ya) {};

\end{tikzpicture}
\caption{The vertex replacement graph $H_7$.}
\label{vertexgraph07}
\end{figure}
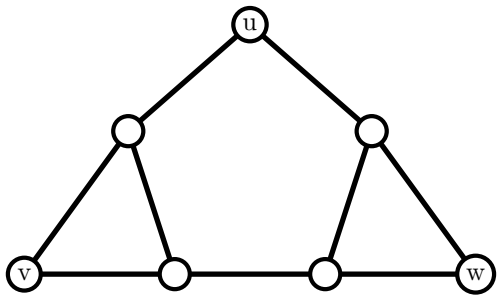

To obtain a graph with no $4$, $8$ or $16$ cycles,
we replace each vertex in $C_{60}$ by a copy
of the graph $H_7$ shown in Figure~\ref{vertexgraph07}.  The three
edges incident with the replaced buckyball vertex are attached to the vertices
labeled $u$, $v$ and $w$ in Figure~\ref{vertexgraph07}.  The attachment
is done so that vertex
$u$ is incident with the double-bond edge.
Since $H_7$ has three vertices of degree $2$, and four vertices of
degree $3$, the resulting graph will evidently be a $3$-connected cubic planar graph
of order $420$.  We denote the graph by $G_{420}$.

It will be convenient to refer to the natural projection from
$V(G_{420})$ to $V(C_{60})$ in which the vertices in each copy of $H_7$ are
projected onto the replaced vertex.

To see that $G_{420}$ contains no $2^m$ cycles for $m \leq 4$, observe
that $H_7$ contains no $2^m$-cycles for any $m$, so any such cycle in $G_{420}$
would project to a cycle in $C_{60}$.
Consider first the cycles
of $G_{420}$ that project to hexagonal face cycles.
Since the minimum distance in $H_7$ between any pair of the attachment
vertices ($u$, $v$ and $w$) is at least $2$, any such cycle
must contain at least $12$ edges from copies of $H_7$, along with
$6$ edges joining different copies, and hence have length at least $18$.
Next consider cycles from $G_{420}$ that project to pentagonal faces.
Since the distance in $H_7$ from $v$ to $w$ is three, each such cycle
contains at least $15$ edges from copies of $H_7$,
along with $5$ edges joining the
copies, and so has length at least $20$.
Finally, because $C_{60}$ contains no other cycles of length less than nine,
one finds that there are no $2^m$ cycles in $G_{420}$ for $m \leq 4$ as claimed.

\section{The general case.}

Next we consider the following variant of \EGC.

\begin{problem}
For $k \ge 3$, what are the smallest cubic graphs with no
$2^m$-cycles, for $m \leq k$?
\end{problem}

Denote the order of a smallest graph by $f(k)$.
It is easy to check that $f(2) = 10$.  There are three
cubic graphs of order $10$ with no $4$-cycles (including the
Petersen graph) and none of smaller order.
Markstr\"{o}m \cite{markstrom} showed that $f(3) = 24$, and listed
all four minimal graphs.
In fact, the graph $H_7$ can viewed as playing a role in one of
Markstr\"{o}m's graphs.
This particular graph can be obtained from $K_4$ by
replacing three of the vertices of $K_4$ by $H_7$ and
replacing the fourth vertex by a copy of $K_3$.

\begin{figure}[H]
\centering
\begin{tikzpicture}[line width=2pt]
\tikzstyle{every node}=[draw=black,fill=white,ultra thick,
  shape=circle,minimum height=0.3cm,inner sep=2];
\def\rad{2.0}
\begin{scope}[rotate=90]
  \foreach \theta in {0,40,...,320}
  {
    \begin{scope}[rotate=\theta]
      \draw (0:3) -- (40:3);
    \end{scope}
  }
  \foreach \theta in {0,120,240}
  {
    \begin{scope}[rotate=\theta]
      \draw (0,0) -- (0:3);
      \draw (-80:3) -- (80:3);
      \draw (0:\rad) -- (0:3);
    \end{scope}
  }
  \foreach \theta in {0,40,...,320}
  {
    \begin{scope}[rotate=\theta]
      \node[fill=white] at (0:3) {};
    \end{scope}
  }
  \node[fill=white] at (0,0) {};
\end{scope}
\end{tikzpicture}
\caption{The Petersen graph.}
\label{petersen3}
\end{figure}
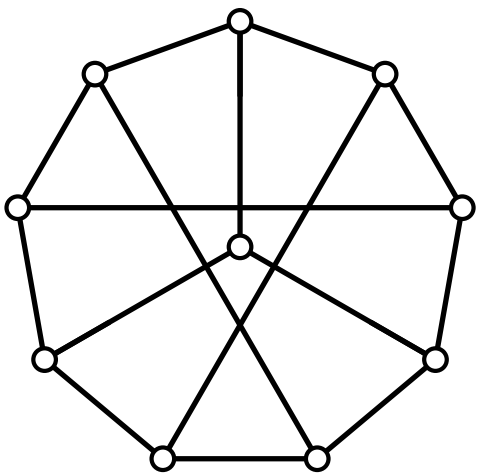

$H_7$ can also be used to construct what appears to be the
smallest known example of a graph with no $2^m$-cycles for
$m \leq 4$.  In this construction, one begins with the Petersen
graph, drawn as in Figure~\ref{petersen3}

\begin{figure}[H]
\centering
\begin{tikzpicture}[line width=2pt]
\tikzstyle{every node}=[draw=black,fill=white,ultra thick,
  shape=circle,minimum height=0.3cm,inner sep=2];
\def\rad{1.5}
\begin{scope}[rotate=90]
  \foreach \theta in {0,40,...,320}
  {
    \begin{scope}[rotate=\theta]
      \draw (0:3) -- (40:3);
    \end{scope}
  }
  \foreach \theta in {0,120,240}
  {
    \begin{scope}[rotate=\theta]
      \draw (0:\rad) -- (120:\rad);
      \draw (-80:3) -- (80:3);
      \draw (0:\rad) -- (0:3);
    \end{scope}
  }
  \foreach \theta in {0,40,...,320}
  {
    \begin{scope}[rotate=\theta]
      \node[fill=white] at (0:3) {};
    \end{scope}
  }
  \foreach \theta in {0,120,240}
  {
    \begin{scope}[rotate=\theta]
      \node[fill=white] at (0:\rad) {};
    \end{scope}
  }
\end{scope}
\end{tikzpicture}
\caption{$G_{12}$, the graph obtained from the Petersen graph by replacing
one vertex by a triangle.}
\label{graph12}
\end{figure}
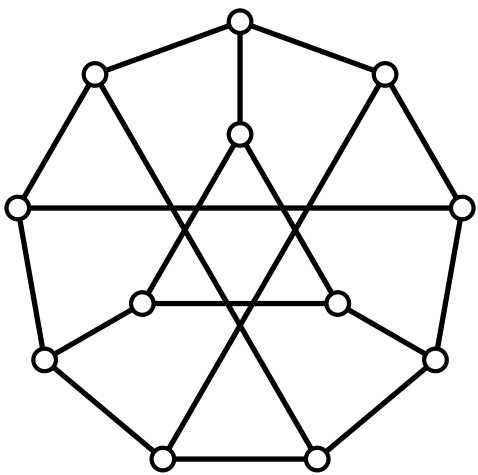

Next we replace the central vertex by a copy of $K_3$ as shown in Figure~\ref{graph12}.

\begin{figure}[H]
\centering
\begin{tikzpicture}
[line width=2,vertex/.style={draw=black,fill=white,ultra thick,
  shape=circle,minimum height=0.3cm,inner sep=2}]
\def\rad{1.5}
\begin{scope}[rotate=90]
  \foreach \theta in {0,40,...,320}
  {
    \begin{scope}[rotate=\theta]
      \draw (0:3) -- (40:3);
    \end{scope}
  }
  
  \draw (  0:3) -- ( 40:3) node[sloped,very near start,fill=black] {};
  \draw ( 40:3) -- (  0:3) node[sloped,very near start,fill=black] {};
  \draw (320:3) -- (  0:3) node[sloped,very near start,fill=black] {};
  \draw (120:3) -- (160:3) node[sloped,very near start,fill=black] {};
  \draw (160:3) -- (120:3) node[sloped,very near start,fill=black] {};
  \draw (200:3) -- (240:3) node[sloped,very near start,fill=black] {};
  \draw (240:3) -- (280:3) node[sloped,very near start,fill=black] {};
  \draw (280:3) -- (240:3) node[sloped,very near start,fill=black] {};
  \draw ( 80:3) -- ( 60:1) node[sloped,very near start,fill=black] {};
  \draw (120:\rad) -- (120:3) node[sloped,very near start,fill=black] {};
  \draw (240:\rad) -- (240:3) node[sloped,very near start,fill=black] {};

  \foreach \theta in {0,120,240}
  {
    \begin{scope}[rotate=\theta]
      \draw (0:\rad) -- (120:\rad);
      \draw (-80:3) -- (80:3);
      \draw (0:\rad) -- (0:3);
    \end{scope}
  }
  \foreach \theta in {0,40,...,320}
  {
    \begin{scope}[rotate=\theta]
      \node[fill=white,vertex] at (0:3) {};
    \end{scope}
  }
  \node[vertex,fill=black] at (0:\rad) {};
  \foreach \theta in {120,240}
  {
    \begin{scope}[rotate=\theta]
      \node[fill=white,vertex] at (0:\rad) {};
    \end{scope}
  }
\end{scope}
\end{tikzpicture}
\caption{Replacing $11$ of the vertices by $H_7$.}
\label{graph78}
\end{figure}
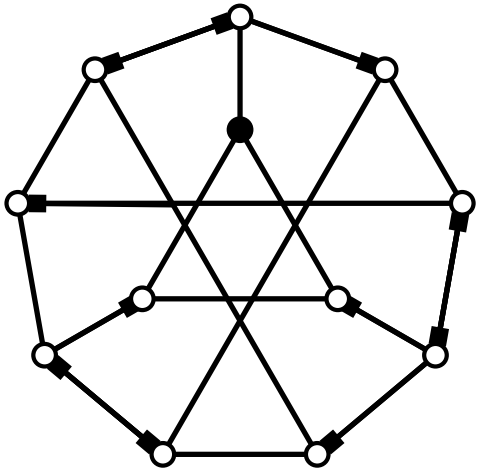

The final step is to replace all but one of the vertices
in $G_{12}$
with a copy of $H_7$.
In Figure~\ref{graph78} we indicate how
this can be done.
The solid vertex in the figure is the vertex
that is {\it not} replaced by a copy of $H_7$.
The heavy boxes attached
to each of the the other vertices
mark the edge that will be
incident with vertex $u$ in that copy of $H_7$.
We thus obtain $G_{78}$, a cubic graph of order $78$.
To see that there are no $4$, $8$ or $16$ cycles in $G_{78}$
one can check all possible cycles that project from $G_{78}$ to 
the single $3$-cycle, the six $5$-cycles and the ten $6$-cycles in
$G_{12}$.  The details are left to the reader.

\begin{figure}[H]
\centering
\begin{tikzpicture}[line width=2pt]
\tikzstyle{every node}=[draw=black,fill=white,ultra thick,
  shape=circle,minimum height=0.3cm,inner sep=2];
\def\rad{2.0}
\begin{scope}[rotate=90]
  \foreach \theta in {6,18,...,354}
  {
    \begin{scope}[rotate=\theta]
      \draw (0:4) -- (12:4);
    \end{scope}
  }
  \foreach \theta in {0,72,...,288}
  {
    \begin{scope}[rotate=\theta]
      \draw ( -6:4) -- (150:4);
      \draw (-18:4) -- ( 90:4);
      \draw (-42:4) -- ( 42:4);
    \end{scope}
  }
  \foreach \theta in {6,18,...,354}
  {
    \begin{scope}[rotate=\theta]
      \node[fill=white] at (0:4) {};
    \end{scope}
  }
\end{scope}
\end{tikzpicture}
\caption{The Tutte-Coxeter graph.}
\label{tuttecoxeter}
\end{figure}
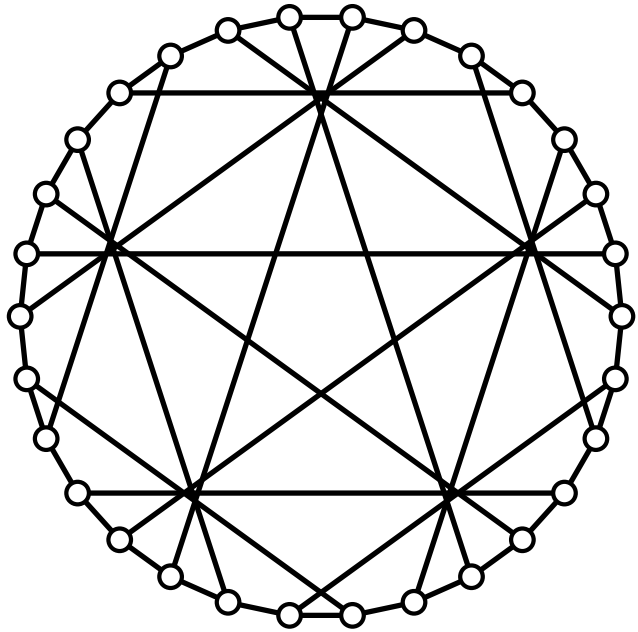

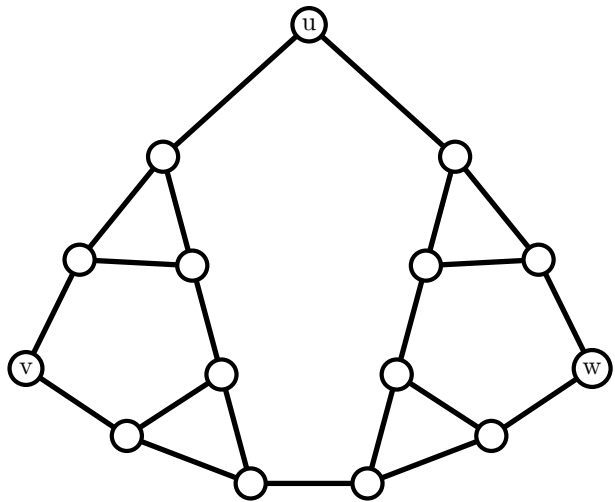
\begin{figure}[H]
\centering
\begin{tikzpicture}[line width=2pt,scale=1.5]
\tikzstyle{every node}=[draw=black,fill=white,ultra thick,
  shape=circle,minimum height=0.4cm,inner sep=2];
\def\xa{{3.0 - cos(72)}}
\def\xb{{4.0 + cos(72)}}
\def\ya{{sin(72)}}
\def\yb{{(cos(36) + 1.0/(2.0 * cos(54)))}}
\def\tha{255}
\def\thb{285}

\begin{scope}[rotate=180]
  \draw (\tha:2) -- (\thb:2);
  \draw (\tha:5) -- (0,-6) -- (\thb:5);

  \begin{scope}[rotate=\tha]
  \draw (\xa,-\ya) -- (2,0) -- (5,0) -- (\xb,-\ya);
  \draw (3,0) -- (\xa,-\ya) -- (3.5,-\yb) -- (\xb,-\ya) -- (4,0);

  \node[fill=white] at (2,0) {};
  \node[fill=white] at (3,0) {};
  \node[fill=white] at (4,0) {};
  \node[fill=white] at (5,0) {};
  \node[fill=white] at (\xa,-\ya) {};
  \node[fill=white] at (3.5,-\yb) {w};
  \node[fill=white] at (\xb,-\ya) {};
  \end{scope}

  \begin{scope}[rotate=\thb]
  \draw (\xa,\ya) -- (2,0) -- (5,0) -- (\xb,\ya);
  \draw (3,0) -- (\xa,\ya) -- (3.5,\yb) -- (\xb,\ya) -- (4,0);

  \node[fill=white] at (2,0) {};
  \node[fill=white] at (3,0) {};
  \node[fill=white] at (4,0) {};
  \node[fill=white] at (5,0) {};
  \node[fill=white] at (\xa,\ya) {};
  \node[fill=white] at (3.5,\yb) {v};
  \node[fill=white] at (\xb,\ya) {};
  \end{scope}

  \node[fill=white] at (0,-6) {u};
\end{scope}

\end{tikzpicture}
\caption{A larger vertex replacement graph $H_{15}$.}
\label{vertexgraph15}
\end{figure}

The final construction gives
a graph with no $2^m$-cycles, for $m \leq 5$.
It is based on the girth $8$ Tutte-Coxeter graph shown in Figure~\ref{tuttecoxeter}
and on the vertex replacement graph $H_{15}$ shown in Figure~\ref{vertexgraph15}.
Note that $H_{15}$ consists of two copies of $H_7$ with one extra vertex, and
that there are no $2^m$-cycles in $H_{15}$.  Note also that the distance
in $H_{15}$ from $u$ to
either $v$ or $w$ is $3$, while the distance from
$v$ to $w$ is $5$.

If one replaces each vertex of Tutte-Coxeter by a copy of $H_{15}$,
there are clearly no $2^m$-cycles for $m \leq 4$, but if the replacement
is done in an arbitary manner, $32$-cycles may result.  To avoid this
one must be a little careful.  Each vertex of Tutte-Coxeter is incident
with two edges on the outer Hamiltonian cycle (in Figure~\ref{tuttecoxeter})
and with one chord edge.  If we replace each Tutte-Coxeter vertex
so that vertex $u$ in
each copy of $H_{15}$ is
incident with a chord edge, then $32$-cycles will be avoided.
This follows
from the observation that
any $8$-cycle in Tutte-Coxeter contains (at least) two consecutive edges on the
outer Hamiltonian cycle, and therefore at least one $v$-$w$ path in a copy
of $H_{15}$.
The resulting graph has order $450$ and no $2^m$-cycles for
$m \leq 5$, and therefore $f(5) \leq 450$.

We summarize the known values and bounds on $f$ in the table at the end
of this note.  The lower bound for $f(4)$ is an unpublished result
of Markstr\"{o}m. 

\begin{center}
\begin{tabular}{r|r} \hline
$k$ & $f(k)$ \\ \hline
$2$ & $10$ \\
$3$ & $24$ \\
$4$ & $54$ -- $78$ \\
$5$ & $ \leq 450$ \\ \hline
\end{tabular}
\end{center}

\end{document}